\documentclass{article}
\usepackage[margin=1in]{geometry}
\usepackage{enumitem}
\usepackage{mathtools}
\usepackage{amssymb}
\usepackage{tikz}
\usetikzlibrary{arrows.meta,graphs,positioning,quotes}
\usepackage[caption=false,font=footnotesize,labelformat=simple]{subfig}
\usepackage{amsthm}
\newtheorem{theorem}{Theorem}
\newtheorem{lemma}[theorem]{Lemma}
\newtheorem{corollary}[theorem]{Corollary}
\newtheorem{proposition}[theorem]{Proposition}
\theoremstyle{definition}
\newtheorem{definition}{Definition}
\newtheorem{example}{Example}

\newtheorem{remark}{Remark}
\usepackage[colorlinks=true,breaklinks=true,bookmarks=true,urlcolor=blue,
     citecolor=blue,linkcolor=blue,bookmarksopen=false,draft=false]{hyperref}

\newenvironment{myproof}[1][Proof]{\begin{proof}[#1.]}{\end{proof}}
\newcommand{\RR}{\mathbb{R}}
\newcommand{\ZZ}{\mathbb{Z}}

\DeclarePairedDelimiter{\norm}{\lVert}{\rVert}

\DeclareMathOperator{\conv}{conv}
\DeclareMathOperator{\dom}{dom}
\DeclareMathOperator{\cl}{cl}

\DeclareMathOperator{\epi}{epi}
\DeclareMathOperator{\st}{s.t.}

\title{Refined Shapley-Folkman Lemma and Its Application in Duality Gap Estimation}
\author{Yingjie Bi and Ao Tang \\
Cornell University, Ithaca, NY, 14850}
\date{}
\begin{document}
\maketitle
\begin{abstract}%
Based on concepts like $k$th convex hull and finer characterization of nonconvexity of a function, we propose a refinement of the Shapley-Folkman lemma and derive a new estimate for the duality gap of nonconvex optimization problems with separable objective functions. We apply our result to a network flow problem and the dynamic spectrum management problem in communication systems as examples to demonstrate that the new bound can be qualitatively tighter than the existing ones. The idea is also applicable to cases with general nonconvex constraints.
\end{abstract}

\section{Introduction}

The Shapley-Folkman lemma (Theorem~\ref{thm:shapleyfolkman}) was stated and used to establish the existence of approximate equilibria in economy with nonconvex preferences \cite{Starr1969}. It roughly says that the sum of a large number of sets is close to a convex set and thus can be used to generalize results on convex objects to nonconvex ones.

\begin{theorem}\label{thm:shapleyfolkman}
Let $S_1,S_2,\dots,S_n$ be any subsets of $\RR^m$. For each $z \in \conv\sum_{i=1}^nS_i=\sum_{i=1}^n\conv S_i$, there exist points $z^i \in \conv S_i$ such that $z=\sum_{i=1}^nz^i$ and $z^i \in S_i$ except for at most $m$ values of $i$.
\end{theorem}

\begin{remark}
In this paper, we use superscript to index vectors and use subscript to refer to a particular component of a vector. For instance, both $x^i$ and $x^{ij}$ are vectors, but $x^i_s$ is the $s$th component of the vector $x^i$. For two vectors $x$ and $y$, $x \leq y$ means $x_s \leq y_s$ holds for all components.
\end{remark}

The Shapley-Folkman lemma has found applications in many fields including economics and optimization theory. In particular, it has been used to study optimization problems with separable objectives and linear constraints such as
\begin{equation}\label{eq:primal}
\begin{split}
\min \quad & \sum_{i=1}^nf_i(x^i) \\
\st \quad & \sum_{i=1}^nA_ix^i \leq b.
\end{split}
\end{equation}
Here $x^i \in \RR^{n_i}$ are the decision variables. The function $f_i:\RR^{n_i} \to \bar\RR$ is proper\footnote{This means that the function never takes $-\infty$ and its domain is nonempty.} and lower semi-continuous, and its domain is bounded. $A_i$ is a matrix of size $m \times n_i$, so there are $m$ linear constraints in total. The dual problem of \eqref{eq:primal} is
\begin{equation}\label{eq:dual}
\begin{split}
\max \quad & -\sum_{i=1}^nf^*_i(-A^T_iy)-b^Ty \\
\st \quad & y \geq 0,
\end{split}
\end{equation}
where $f^*_i$ is the conjugate function of $f_i$. Denote the optimal value of the primal problem \eqref{eq:primal} and dual problem \eqref{eq:dual} as $p$ and $d$, respectively. In general, there will be a positive duality gap $p-d>0$ if some function $f_i$ is not convex.

The authors of \cite{AE1976} presented the following upper bound for the duality gap based on the Shapley-Folkman lemma:
\begin{equation}\label{eq:originalbound}
p-d \leq \min\{m+1,n\}\max_{i=1,\dots,n}\rho(f_i).
\end{equation}
Here $\rho(f)$ is the nonconvexity of a proper function $f$ defined by
\begin{equation}\label{eq:nonconvexity}
\rho(f)=\sup\left\{f\left(\sum_j\alpha_jx^j\right)-\sum_j\alpha_jf(x^j)\right\}
\end{equation}
over all finite convex combinations of points $x^j \in \dom f$, i.e., $f(x^j)<+\infty$, $\alpha_j \geq 0$ with $\sum_j\alpha_j=1$.

In \cite{UB2016}, an improved bound for the duality gap\footnote{In fact, the bound \eqref{eq:udellbound} derived in \cite{UB2016} is for the difference $p-\hat p$, in which $\hat p$ is the optimal value of the convexified problem where the function $f_i$ is substituted by $f^{**}_i$. However, $\hat p=d$ because the refined Slater's condition holds for the convexified problem.} was given by
\begin{equation}\label{eq:udellbound}
p-d \leq \sum_{i=1}^{\min\{m,n\}}\rho(f_i),
\end{equation}
where we assume that $\rho(f_1) \geq \dots \geq \rho(f_n)$. Although the new bound \eqref{eq:udellbound} is only a slight improvement over the original bound \eqref{eq:originalbound} by a factor of $m/(m+1)$, it nevertheless demonstrates that \eqref{eq:originalbound} can never be tight except for some trivial situations. Thus, using Shapley-Folkman lemma as in \cite{AE1976} does not result in a tight bound for the duality gap, which suggests that the original Shapley-Folkman lemma itself can be improved in certain circumstances.

In this paper, we propose a refinement for the original Shapley-Folkman lemma and derive a new bound for the duality gap based on this new result. The refined Shapley-Folkman lemma is stated and proved in Section~\ref{sec:refined}. Unlike \eqref{eq:originalbound} and \eqref{eq:udellbound}, our new bound for the duality gap depends on some finer characterization of the nonconvexity of a function, which is introduced in Section~\ref{sec:nonconvexity}. The new bound itself is given in Section~\ref{sec:dualitygap}, which can easily recover existing ones like \eqref{eq:udellbound}. Next, we will apply it to two examples, a network flow problem and the dynamic spectrum management problem, in Section~\ref{sec:application} to demonstrate that the new bound can be qualitatively tighter than the bound \eqref{eq:udellbound}. Although we mainly focus on the case of linear constraints, Section~\ref{sec:nonlinear} shows how the major idea in this paper can be applied to the cases with general convex or even nonconvex constraints.

\section{Refined Shapley-Folkman Lemma}\label{sec:refined}

To write down our refined version of the Shapley-Folkman lemma, we need to first introduce the concept of $k$th convex hull.

\begin{definition}
The \emph{$k$th convex hull} of a set $S$, denoted by $\conv_k S$, is the set of convex combinations of $k$ points in $S$, i.e.,
\[
\conv_k S=\left\{\sum_{j=1}^k\alpha_jv^j \middle| v^j \in S, \alpha_j \geq 0, \forall j=1,\dots,k, \sum_{j=1}^k\alpha_j=1\right\}.
\]
\end{definition}

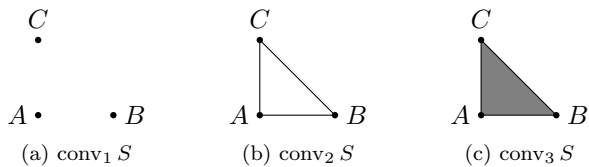
\begin{figure}
\centering
\tikzset{every node/.style={circle,fill=black,inner sep=0.03cm}}
\subfloat[$\conv_1 S$]{
\begin{tikzpicture}
\draw (0,0) node["$A$" left]{} (1,0) node["$B$" right]{} (0,1) node["$C$" above]{};
\end{tikzpicture}
}\qquad
\subfloat[$\conv_2 S$]{
\begin{tikzpicture}
\draw (0,0) node["$A$" left]{} -- (1,0) node["$B$" right]{} -- (0,1) node["$C$" above]{} -- cycle;
\end{tikzpicture}\label{fig:conv2s}
}\qquad
\subfloat[$\conv_3 S$]{
\begin{tikzpicture}
\draw[fill=gray] (0,0) node["$A$" left]{} -- (1,0) node["$B$" right]{} -- (0,1) node["$C$" above]{} -- cycle;
\end{tikzpicture}\label{fig:conv3s}
}
\caption{The $k$th convex hull of a three-point set $S=\{A,B,C\}$.}\label{fig:convks}
\end{figure}

Figure~\ref{fig:convks} gives a simple example to illustrate the definition of $k$th convex hull. In Figure~\ref{fig:convks}, the set $S=\{A,B,C\}$, $\conv_1 S=S$, $\conv_2 S$ are the segments $AB$, $BC$ and $CA$, while $\conv_3 S$ is the full triangle which is also the convex hull of set $S$. In general, Carath\'eodory's theorem implies that $\conv_{m+1}S=\conv S$ for any set $S \subseteq \RR^m$. However, for a particular set, the minimum $k$ such that $\conv_k S=\conv S$ can be smaller than $m+1$, and this number intuitively reflects how the set is closer to being convex. For instance, if we start from $T=\conv_2 S$, the set in Figure~\ref{fig:conv2s}, then $\conv_k T=\conv T$ for $k=2$.

Next, we recall the concept of $k$-extreme points of a convex set, which is a generalization of extreme points.

\begin{definition}
A point $z$ in a convex set $S$ is called a \emph{$k$-extreme} point of $S$ if we cannot find $(k+1)$ independent vectors $d^1,d^2,\dots,d^{k+1}$ such that $z \pm d^i \in S$.
\end{definition}

According to our definition, if a point is $k$-extreme, then it is also $k'$-extreme for $k' \geq k$. For a convex set in $\RR^m$, a point is an extreme point if and only if it is 0-extreme, a point is on the boundary if and only if it is $(m-1)$-extreme, and every point is $m$-extreme. For example, in Figure~\ref{fig:conv3s}, the vertices $A,B,C$ are 0-extreme, the points on segments $AB$, $BC$ and $CA$ are 1-extreme, and all the points are 2-extreme.

Now we can state our refined Shapley-Folkman lemma:

\begin{theorem}\label{thm:refined}
Let $S_1,S_2,\dots,S_n$ be any subsets of $\RR^m$. Assume $z$ is a $k$-extreme point of $\conv\sum_{i=1}^nS_i$, then there exist integers $1 \leq k_i \leq k+1$ with $\sum_{i=1}^nk_i \leq k+n$ and points $z^i \in \conv_{k_i}S_i$ such that $z=\sum_{i=1}^nz^i$.
\end{theorem}

The original Shapley-Folkman lemma (Theorem~\ref{thm:shapleyfolkman}) now becomes a direct corollary of Theorem~\ref{thm:refined}. Since any point $z \in \conv\sum_{i=1}^nS_i$ is an $m$-extreme point. Applying Theorem~\ref{thm:refined} on this point gives a decomposition $z=\sum_{i=1}^nz^i$ with $z^i \in \conv_{k_i}S_i \subseteq \conv S_i$ and $\sum_{i=1}^nk_i \leq m+n$. Then the conclusion in Theorem~\ref{thm:shapleyfolkman} follows because $z^i \in S_i$ if $k_i=1$, while the number of indices $i$ with $k_i \geq 2$ is bounded by $m$.

\begin{remark}
Our Theorem~\ref{thm:refined} is similar to the refined version of Shapley-Folkman lemma proposed in \cite{LS2013}. However, the result in \cite{LS2013} does not take extremeness of the point into account. By Theorem~\ref{thm:refined}, stronger result can be obtained from the knowledge of extremeness.
\end{remark}

To prove Theorem~\ref{thm:refined}, we need the following property of $k$-extreme points in a polyhedron:

\begin{lemma}\label{lem:kexposed}
Let $P \subseteq \RR^m$ be a polyhedron and $z$ be a $k$-extreme point of $P$, then there exists a vector $a \in \RR^m$ such that the set $\{y \in P | a^Ty \leq a^Tz\}$ is in a $k$-dimensional affine subspace.
\end{lemma}
\begin{myproof}
Assume that the polyhedron $P$ is represented by $Ax \geq b$. Let $A_=$ be the submatrix of $A$ containing the rows of active constraints for the point $z$, and $b_=$ be the vector containing the corresponding constants in $b$. The dimension of the kernel of $A_=$ is at most $k$. Otherwise, we can find independent and sufficiently small vectors $d^1,\dots,d^{k+1}$ such that $A_=d^i=0$ and $A(z \pm d^i) \geq b$ for $i=1,\dots,k+1$. This implies $z \pm d^i \in P$, which contradicts with the $k$-extremeness of point $z$.

Let $a$ be the vector such that $a^T$ is the sum of all rows in $A_=$. Consider a point $y$ satisfying $Ay \geq b$ and $a^Ty \leq a^Tz$. Since adding all inequalities together in $A_=y \geq b_==A_=z$ gives $a^Ty \geq a^Tz$, we must have $A_=y=b_=$. Therefore, $y$ is in the affine subspace defined by $A_=x=b_=$ whose dimension is at most $k$.
\end{myproof}

\begin{remark}
In the literature, the point satisfying the conclusion of Lemma~\ref{lem:kexposed} is called a \emph{$k$-exposed point}. For a general convex set $S$, a $k$-extreme point may fail to be a $k$-exposed point, although it must be in the closure of the set of $k$-exposed points if $S$ is compact \cite{Asplund1963}. For the special case of polyhedra, these two concepts are equivalent, and Lemma~\ref{lem:kexposed} is a generalization of the well-known result that an extreme point of a polyhedron is the unique minimizer of some linear function.
\end{remark}

\begin{myproof}[Proof of Theorem~\ref{thm:refined}]
Since $z$ is in the convex hull of $\sum_{i=1}^nS_i$, there exists some integer $l$ such that $z$ can be written as
\begin{equation}\label{eq:zinconv}
z=\sum_{j=1}^l\alpha_j\sum_{i=1}^nv^{ij},
\end{equation}
in which $v^{ij} \in S_i$, $\alpha_j \geq 0$, $j=1,\dots,l$ and $\sum_{j=1}^l\alpha_j=1$.

Define $S'_i=\{v^{i1},\dots,v^{il}\} \subseteq S_i$, then \eqref{eq:zinconv} actually tells us that $z \in \conv\sum_{i=1}^nS'_i$, so $z$ must be $k$-extreme in this polytope that lies in $\conv\sum_{i=1}^nS_i$. By Lemma~\ref{lem:kexposed}, there exists a vector $a \in \RR^m$ such that the set
\[
\left\{y \in \conv\sum_{i=1}^nS'_i \middle| a^Ty \leq a^Tz\right\}
\]
is in a $k$-dimensional affine subspace $L$ of $\RR^m$. Without loss of generality, we assume that the subspace
\[
L=\{y \in \RR^m | y_{k+1}=y_{k+2}=\dots=y_m=0\}.
\]

Next, consider the following linear program in which $\beta_{ij}$ are the decision variables:
\[
\begin{split}
\min \quad & \sum_{i=1}^n\sum_{j=1}^l\beta_{ij}a^Tv^{ij} \\
\st \quad & \sum_{i=1}^n\sum_{j=1}^l\beta_{ij}v^{ij}_s=z_s, \quad \forall s=1,\dots,k, \\
& \sum_{j=1}^l\beta_{ij}=1, \quad \forall i=1,\dots,n, \\
& \beta_{ij} \geq 0, \quad \forall i=1,\dots,n, \: \forall j=1,\dots,l.
\end{split}
\]
Setting $\beta_{ij}=\alpha_j$ gives a feasible solution to the above problem with objective value $a^Tz$. Among all the optimal solutions, pick up a particular vertex solution $\beta^*_{ij}$, which should have at least $nl$ active constraints. We already have $k+n$ active constraints, so the number of nonzero $\beta^*_{ij}$ entries is at most $k+n$. Define
\[
z^i=\sum_{j=1}^l\beta^*_{ij}v^{ij}, \quad z'=\sum_{i=1}^nz^i,
\]
and let $k_i$ be the number of nonzero entries in $\beta^*_{i1},\dots,\beta^*_{il}$. Since $\sum_{j=1}^l\beta^*_{ij}=1$, there must be a nonzero one and thus $k_i \geq 1$. Now we know that $z^i \in \conv_{k_i}S_i$, and $\sum_{i=1}^nk_i \leq k+n$ implies that each $k_i$ cannot exceed $k+1$. The remaining thing to show is $z_s=z'_s$ for $s=k+1,\dots,m$. Because
\[
z' \in \sum_{i=1}^n\conv S'_i=\conv\sum_{i=1}^nS'_i
\]
and
\[
a^Tz'=\sum_{i=1}^n\sum_{j=1}^l\beta^*_{ij}a^Tv^{ij} \leq a^Tz,
\]
$z' \in L$. Since $z \in L$, the last $m-k$ components of both $z$ and $z'$ are all zeros, so $z=z'$.
\end{myproof}

A useful consequence can be made from Theorem~\ref{thm:refined}:

\begin{corollary}\label{cor:refinedhalf}
Let $S_1,S_2,\dots,S_n$ be any subsets of $\RR^m$. If $z \in \conv\sum_{i=1}^nS_i$, then there exist integers $1 \leq k_i \leq m$ with $\sum_{i=1}^nk_i \leq m-1+n$ and points $z^i \in \conv_{k_i}S_i$ such that $z_s=\sum_{i=1}^nz^i_s$ for $s=1,\dots,m-1$ and $z_m \geq \sum_{i=1}^nz^i_m$.
\end{corollary}
\begin{myproof}
Using the same argument in the proof of Theorem~\ref{thm:refined}, choose $S'_i \subseteq S_i$ containing finite points such that $z \in \conv\sum_{i=1}^nS'_i$. Since $\conv\sum_{i=1}^nS'_i$ is a compact set,
\[
\inf\left\{w_m \middle| w \in \conv\sum_{i=1}^nS'_i, w_1=z_1, \dots, w_{m-1}=z_{m-1}\right\}
\]
can be achieved by some point $w^*$. $w^*$ is an $(m-1)$-extreme point of $\conv\sum_{i=1}^nS'_i$, and applying Theorem~\ref{thm:refined} on the point $w^*$ gives the desired result.
\end{myproof}

In Section~\ref{sec:dualitygap}, when we prove the bound for the duality gap, Corollary~\ref{cor:refinedhalf} is applied to the epigraph of an $m$-variable function, where $m$ is the number of constraints. The bound in Corollary~\ref{cor:refinedhalf} will be $m+n$ instead of $m+1+n$ from Theorem~\ref{thm:refined} directly. Although the last component in the decomposition given by Corollary~\ref{cor:refinedhalf} is only an inequality, we will see in the following that the inequality is sufficient for our purpose. Therefore, Corollary~\ref{cor:refinedhalf} is the true source of the improvement in \cite{UB2016} from the original bound \eqref{eq:originalbound}.

\section{Characterization of Nonconvexity}\label{sec:nonconvexity}

To improve the bound \eqref{eq:udellbound}, some finer characterization of the nonconvexity of a function has to be introduced. In parallel to the definition of $k$th convex hull of a set, define the \emph{$k$th nonconvexity} $\rho^k(f)$ of a proper function $f$ to be the supremum in \eqref{eq:nonconvexity} taken over the convex combinations of $k$ points instead of arbitrary number of points. Obviously,
\[
0=\rho^1(f) \leq \rho^2(f) \leq \dots \leq \rho(f).
\]
In fact, we have the following property:

\begin{proposition}\label{prop:nonconvexity}
For any proper function $f:\RR^n \to \bar\RR$, $\rho^{n+1}(f)=\rho(f)$.
\end{proposition}
\begin{myproof}
We only need to show that $\rho(f) \leq \rho^{n+1}(f)$. Choose any convex combination $x=\sum_{j=1}^l \alpha_jx^j$ with all points $x^j \in \dom f$, $\alpha_j \geq 0$, and $\sum_{j=1}^l \alpha_j=1$. Since $(x^j,f(x^j)) \in \epi f$, the point
\[
\left(\sum_{j=1}^l \alpha_jx^j,\sum_{j=1}^l \alpha_jf(x^j)\right) \in \conv\epi f.
\]
Using Corollary~\ref{cor:refinedhalf}, we can find $(y^i,t_i) \in \epi f$, $\beta_i \geq 0$, $i=1,\dots,n+1$, and $\sum_{i=1}^{n+1} \beta_i=1$ such that
\[
x=\sum_{j=1}^l \alpha_jx^j=\sum_{i=1}^{n+1} \beta_iy^i, \quad \sum_{j=1}^l \alpha_jf(x^j) \geq \sum_{i=1}^{n+1} \beta_it_i.
\]
Now
\[
f(x)-\sum_{j=1}^l \alpha_jf(x^j) \leq f\left(\sum_{i=1}^{n+1} \beta_iy^i\right)-\sum_{i=1}^{n+1}\beta_it_i \leq f\left(\sum_{i=1}^{n+1} \beta_iy^i\right)-\sum_{i=1}^{n+1}\beta_if(y^i),
\]
which implies $\rho(f) \leq \rho^{n+1}(f)$.
\end{myproof}

For lower semi-continuous functions, the following proposition provides an equivalent definition for the $k$th nonconvexity, which sheds light on the connection between the concepts of $k$th nonconvexity and $k$th convex hull.

\begin{proposition}\label{prop:nonconvexityalt}
Assume a proper function $f$ is lower semi-continuous and bounded below by some affine function. Let $f^{(k)}$ be the function whose epigraph is the closure of the $k$th convex hull of the epigraph of $f$, i.e.,
\[
\epi f^{(k)}=\cl\conv_k\epi f.
\]
Then
\begin{equation}\label{eq:nonconvexityalt}
\rho^k(f)=\sup_x\{f(x)-f^{(k)}(x)\},
\end{equation}
where we interpret $(+\infty)-(+\infty)=0$.
\end{proposition}
\begin{myproof}
The assumption on the function $f$ implies that $f^{(k)}$ is also a proper function. Consider an arbitrary $k$-point convex combination of points $x^j \in \dom f$, for $j=1,\dots,k$. Following the first step in the proof of Proposition~\ref{prop:nonconvexity}, we have
\[
\left(\sum_{j=1}^k \alpha_jx^j,\sum_{j=1}^k \alpha_jf(x^j)\right) \in \conv_k\epi f \subseteq \cl\conv_k\epi f=\epi f^{(k)}.
\]
Therefore,
\[
f\left(\sum_{j=1}^k \alpha_jx^j\right)-\sum_{j=1}^k \alpha_jf(x^j) \leq f\left(\sum_{j=1}^k \alpha_jx^j\right)-f^{(k)}\left(\sum_{j=1}^k \alpha_jx^j\right),
\]
which implies
\[
\rho^k(f) \leq \sup_x\{f(x)-f^{(k)}(x)\}.
\]

On the other hand, for any $x \in \dom f^{(k)}$,
\[
(x,f^{(k)}(x)) \in \epi f^{(k)}=\cl\conv_k\epi f.
\]
In the case of $x \in \dom f$, by the lower semi-continuity of $f$, for every $\epsilon>0$, there exists $(\kappa,\eta) \in \conv_k\epi f$ which is sufficiently close to $(x,f^{(k)}(x))$ such that
\[
f(\kappa) \geq f(x)-\epsilon, \quad \eta \leq f^{(k)}(x)+\epsilon.
\]
Because $(\kappa,\eta) \in \conv_k\epi f$, there exists $\alpha_j \geq 0$ for $j=1,\dots,k$ such that $\sum_{j=1}^k\alpha_j=1$ and
\[
\kappa=\sum_{j=1}^k\alpha_jx^j, \quad \eta \geq \sum_{j=1}^k\alpha_jf(x^j)
\]
in which $x^j \in \dom f$. Thus
\[
f(x)-f^{(k)}(x) \leq f(\kappa)-\eta+2\epsilon \leq f\left(\sum_{j=1}^k\alpha_jx^j\right)-\sum_{j=1}^k\alpha_jf(x^j)+2\epsilon.
\]
The $f(x)=+\infty$ case can be dealt with similarly, and we can conclude
\[
\sup_x\{f(x)-f^{(k)}(x)\} \leq \rho^k(f)
\]
by letting $\epsilon \to 0$.
\end{myproof}

\begin{remark}
If a proper function $f$ is bounded below by some affine function, then
\[
\epi f^{**}=\cl\conv f
\]
(see \cite[Theorem X.1.3.5]{HL1993}). Therefore, \eqref{eq:nonconvexityalt} can be regarded as a generalization for the alternative definition of nonconvexity
\[
\rho(f)=\sup_x\{f(x)-f^{**}(x)\}
\]
used in \cite{UB2016}.
\end{remark}

In the remaining of this section, three examples will be given to illustrate how to calculate the $k$th nonconvexity of a particular function. The results will be used in Section~\ref{sec:application}.

\begin{example}\label{exmp:linear}
Consider the function
\[
f(x)=f(x_1,\dots,x_n)=\min_{s=1,\dots,n}x_s
\]
defined on the box $0 \leq x \leq 1$, $x \in \RR^n$. It is already known that $\rho(f)=(n-1)/n$ (see \cite[Table 1]{UB2016}). By Proposition~\ref{prop:nonconvexity}, $\rho^k(f)=\rho(f)=(n-1)/n$ for $k \geq n+1$.

For $k=1,\dots,n$, as in the proof of Proposition~\ref{prop:nonconvexityalt}, pick up any $k$-point convex combination of points $0 \leq x^j \leq 1$, $j=1,\dots,k$. For a given $i \in \{1,\dots,k\}$, let $s(i)$ be the index such that $x^i_{s(i)}$ is the minimum among $x^i_1,\dots,x^i_n$, then
\[
\begin{split}
f(x)&=\min_{s=1,\dots,n}\left\{\sum_{j=1}^k \alpha_jx^j_s\right\} \leq \sum_{j=1}^k \alpha_jx^j_{s(i)} \\
    &\leq \alpha_ix^i_{s(i)}+1-\alpha^i=\alpha_if(x^i)+1-\alpha^i.
\end{split}
\]
Summing up among $i=1,\dots,k$, we have
\[
kf(x) \leq \sum_{i=1}^k \alpha_if(x^i)+k-1,
\]
which implies
\[
f(x)-\sum_{i=1}^k \alpha_if(x^i) \leq \frac{k-1}{k}\left(1-\sum_{i=1}^k \alpha_if(x^i)\right) \leq \frac{k-1}{k}.
\]
The above argument shows that $\rho^k(f) \leq (k-1)/k$. In fact, the equality holds, which can be easily seen by considering the average of first $k$ points of
\begin{gather*}
x^1=(0,1,\dots,1), \\
x^2=(1,0,\dots,1), \\
\dots, \\
x^n=(1,1,\dots,0).
\end{gather*}
In conclusion,
\[
\rho^k(f)=\begin{dcases*}
\frac{k-1}{k}, & if $k=1,\dots,n$, \\
\frac{n-1}{n}, & if $k \geq n+1$.
\end{dcases*}
\]
\end{example}

\begin{example}\label{exmp:log}
Consider the function
\[
g(x)=g(x_1,\dots,x_n)=-\log\max_{s=1,\dots,n}x_s
\]
defined on the region $x \geq 0$ except $x=0$.

For $k=1,\dots,n$, pick up any $k$-point convex combination. Without loss of generality, assume the coefficients $\alpha_j>0$ for $j=1,\dots,k$. For a given $i \in \{1,\dots,k\}$, let $s(i)$ be the index such that $x^i_{s(i)}$ is the maximum among $x^i_1,\dots,x^i_n$, then
\[
\begin{split}
g(x)&=-\log\max_{s=1,\dots,n}\left\{\sum_{j=1}^k \alpha_jx^j_s\right\} \leq -\log\sum_{j=1}^k \alpha_jx^j_{s(i)} \\
    &\leq -\log(\alpha_ix^i_{s(i)})=-\log\alpha_i+g(x^i).
\end{split}
\]
Summing up among $i=1,\dots,k$ with weight $\alpha_i$, we have
\[
\begin{split}
g(x) &\leq -\sum_{i=1}^k \alpha_i\log\alpha_i+\sum_{i=1}^k \alpha_ig(x^i) \\
     &\leq \log k+\sum_{i=1}^k \alpha_ig(x^i).
\end{split}
\]
The above argument shows that $\rho^k(g) \leq \log k$. In fact, the equality holds, which can be easily seen by considering the average of first $k$ points of
\begin{gather*}
x^1=(1,0,\dots,0), \\
x^2=(0,1,\dots,0), \\
\dots, \\
x^n=(0,0,\dots,1).
\end{gather*}

To calculate $\rho^{n+1}(g)$, define $h(x)=-\log\sum_{s=1}^n x_s$. Then $h(x)$ is convex and $g(x)-\log n \leq h(x) \leq g(x)$. Thus, for any $(n+1)$-point convex combination,
\[
\begin{split}
g\left(\sum_{j=1}^{n+1} \alpha_jx^j\right) &\leq h\left(\sum_{j=1}^{n+1} \alpha_jx^j\right)+\log n \\
                                           &\leq \sum_{j=1}^{n+1} \alpha_jh(x^j)+\log n \\
                                           &\leq \sum_{j=1}^{n+1} \alpha_jg(x^j)+\log n.
\end{split}
\]
Therefore, $\rho^{n+1}(g) \leq \log n$. On the other hand, $\rho^{n+1}(g) \geq \rho^n(g)=\log n$.

In conclusion,
\[
\rho^k(g)=\begin{dcases*}
\log k, & if $k=1,\dots,n$, \\
\log n, & if $k \geq n+1$.
\end{dcases*}
\]
\end{example}

For more complex functions, it is usually hard to compute its $k$th nonconvexity exactly. However, sometimes we can approximate the $k$th nonconvexity of a function by reducing it to another function whose nonconvexity is already known. This technique is demonstrated by the following example:

\begin{example}\label{exmp:capacity}
Consider the function
\[
h_\sigma(x)=h_\sigma(x_1,\dots,x_n)=\sum_{s=1}^n\log\frac{\norm{x}_1-x_s+\sigma}{\norm{x}_1+\sigma}
\]
defined on the box $0 \leq x \leq 1$, $x \in \RR^n$. Here $\sigma$ is a parameter in the range $0<\sigma \leq 1$.

Define an auxiliary function
\[
H(x;\sigma)=\prod_{s=1}^n\frac{\norm{x}_1-x_s+\sigma}{\norm{x}_1+\sigma},
\]
then $h_\sigma(x)=\log H(x;\sigma)$. To compute the $k$th nonconvexity for the function $h_\sigma$, we first prove some elementary properties for the function $H(x;\sigma)$.

\begin{lemma}\label{lem:hprop}
The function $H(x;\sigma)$ has the following properties:
\begin{enumerate}[label=(\alph*)]
\item For any vectors $x$ and $y$ in the region $0 \leq x,y \leq \sigma$, if $y \leq x$, then $H(y;\sigma) \geq H(x;\sigma)$.
\item $\sigma H(x;1) \leq H(x;\sigma) \leq H(x;1)$.
\end{enumerate}
\end{lemma}
\begin{myproof}
For any $x$ in the region $0 \leq x \leq \sigma$, the partial derivatives
\begin{align*}
\frac{\partial H(x;\sigma)}{\partial x_i}&=H(x;\sigma)\left(\sum_{s=1}^n\frac{1}{\norm{x}_1-x_s+\sigma}-\frac{1}{\norm{x}_1-x_i+\sigma}-\frac{n}{\norm{x}_1+\sigma}\right) \\
&=H(x;\sigma)\left(\sum_{s=1}^n\frac{x_s}{(\norm{x}_1-x_s+\sigma)(\norm{x}_1+\sigma)}-\frac{1}{\norm{x}_1-x_i+\sigma}\right) \\
&\leq H(x;\sigma)\left(\sum_{s=1}^n\frac{x_s}{\norm{x}_1(\norm{x}_1+\sigma)}-\frac{1}{\norm{x}_1+\sigma}\right)=0,
\end{align*}
which gives the first property.

For the second property, it is obvious to see that $H(x;\sigma) \leq H(x;1)$. The other inequality is equivalent to
\[
p(\sigma)=\frac{1}{\sigma}H(x;\sigma)-H(x;1) \geq 0.
\]
The partial derivative
\begin{align*}
\frac{\partial H(x;\sigma)}{\partial\sigma}&=H(x;\sigma)\left(\sum_{s=1}^n\frac{1}{\norm{x}_1-x_s+\sigma}-\frac{n}{\norm{x}_1+\sigma}\right) \\
&=H(x;\sigma)\sum_{s=1}^n\frac{x_s}{(\norm{x}_1-x_s+\sigma)(\norm{x}_1+\sigma)} \\
&\leq H(x;\sigma)\sum_{s=1}^n\frac{x_s}{\sigma(\norm{x}_1+\sigma)}=\frac{1}{\sigma}H(x;\sigma)\frac{\norm{x}_1}{\norm{x}_1+\sigma} \leq \frac{1}{\sigma}H(x;\sigma)
\end{align*}
implies that
\[
p'(\sigma)=-\frac{1}{\sigma^2}H(x;\sigma)+\frac{1}{\sigma}\frac{\partial H(x;\sigma)}{\partial\sigma} \leq 0.
\]
Therefore, the function $p(\sigma)$ is nonincreasing. Together with $p(1)=0$, we have proved the nonnegativity of $p(\sigma)$.
\end{myproof}

To upper bound the $k$th nonconvexity of the function $h_\sigma$, consider arbitrary points $x^j$ for $j=1,\dots,k$ with corresponding combination weights $\alpha_j>0$. Define $k$ vectors $y^1,\dots,y^k$ in $\RR^k$ by
\begin{gather*}
y^1=(1/H(x^1;1),0,\dots,0), \\
y^2=(0,1/H(x^2;1),\dots,0), \\
\dots, \\
y^k=(0,0,\dots,1/H(x^k;1)).
\end{gather*}
Using the result of nonconvexity for the function $g$ given in Example~\ref{exmp:log} and the properties proved in Lemma~\ref{lem:hprop}, we have
\begin{align*}
h_\sigma&\left(\sum_{j=1}^k\alpha_jx^j\right)=\log H\left(\sum_{j=1}^k\alpha_jx^j;\sigma\right) \leq \log H\left(\sum_{j=1}^k\alpha_jx^j;1\right) && \text{by property (b)} \\
&\leq \log\max_{j=1,\dots,k}H(\alpha_jx^j;1) && \text{by property (a)} \\
&\leq \log\max_{j=1,\dots,k}\frac{1}{\alpha_j}H(\alpha_jx^j;\alpha_j)=\log\max_{j=1,\dots,k}\frac{1}{\alpha_j}H(x^j;1) && \text{by property (b)} \\
&=g\left(\sum_{j=1}^k\alpha_jy^j\right) \leq \sum_{j=1}^k\alpha_jg(y^j)+\log k && \text{by the nonconvexity of $g$} \\
&=\sum_{j=1}^k\alpha_j\log H(x^j;1)+\log k \leq \sum_{j=1}^k\alpha_j\log\frac{1}{\sigma}H(x^j;\sigma)+\log k && \text{by property (b)} \\
&=\sum_{j=1}^k\alpha_jh_\sigma(x^j)+\log\frac{k}{\sigma}.
\end{align*}
The above argument shows that the $k$th nonconvexity $\rho^k(h_\sigma) \leq \log(k/\sigma)$.
\end{example}

\begin{remark}
In Example~\ref{exmp:capacity} above, an upper bound for the $k$th nonconvexity of function $h_\sigma$ is obtained by a reduction from the nonconvexity of $g$ in Example~\ref{exmp:log}. Along this line of thoughts, it is conceivable to find the exact value for the $k$th nonconvexity of $h_\sigma$ if we are able to reduce $h_\sigma$ to itself (but with just $k$ variables).
\end{remark}

\section{Bounding Duality Gap}\label{sec:dualitygap}

Now we can state the main result on the duality gap between the primal problem \eqref{eq:primal} and the dual problem \eqref{eq:dual}.

\begin{theorem}\label{thm:dualitygap}
Assume that the primal problem \eqref{eq:primal} is feasible, i.e., $p<+\infty$. Then there exist integers $1 \leq k_i \leq m+1$ such that $\sum_{i=1}^n k_i \leq m+n$ and the duality gap
\[
p-d \leq \sum_{i=1}^n \rho^{k_i}_i.
\]
Here $\rho^k_i=\rho^k(f_i)$ is the $k$th nonconvexity of function $f_i$.
\end{theorem}

First, let us define the perturbation function $v: \RR^m \to \bar\RR$ by letting $v(z)$ be the optimal value of the perturbed problem
\[
\begin{split}
\min \quad & \sum_{i=1}^n f_i(x^i) \\
\st \quad  & \sum_{i=1}^n A_ix^i \leq b+z.
\end{split}
\]
As is the case with convex optimization, $p=v(0)$ and $d=v^{**}(0)$ (see \cite[Lemma~2.2]{ET1999}).

\begin{lemma}\label{lem:perturbfunc}
The perturbation function $v$ is lower semi-continuous.
\end{lemma}

\begin{myproof}
Pick any $z \in \RR^m$. We want to show that if $z^k \to z$ as $k \to \infty$,
\[
l=\liminf_{k \to \infty}v(z^k) \geq v(z).
\]
The above inequality clearly holds when $l=+\infty$. If $l<+\infty$, by considering a subsequence of $\{v(z^k)\}_{k=1}^\infty$, without loss of generality we can assume $v(z^k)<+\infty$ for each $k$ and
\[
\lim_{k \to \infty}v(z^k)=l.
\]
For each $k$, find $(\hat x^{1k},\dots,\hat x^{nk})$ attaining the optimal value of the perturbed problem related to $v(z^k)$, i.e.,
\[
v(z^k)=\sum_{i=1}^n f_i(\hat x^{ik}), \quad \sum_{i=1}^n A_i\hat x^{ik} \leq b+z^k.
\]
By extracting a convergent subsequence for each $\{\hat x^{ik}\}_{k=1}^\infty$, we can assume $\{\hat x^{ik}\}_{k=1}^\infty$ has a limit $x^i$. Then
\[
\sum_{i=1}^n A_ix^i \leq b+z,
\]
which implies that $(x^1,\dots,x^n)$ is feasible to the perturbed problem related to $v(z)$, so
\[
\sum_{i=1}^n f_i(x^i) \geq v(z).
\]
Now
\[
l=\lim_{k \to \infty}\sum_{i=1}^n f_i(\hat x^{ik}) \geq \sum_{i=1}^n \liminf_{k \to \infty}f_i(\hat x^{ik}) \geq \sum_{i=1}^n f_i(x^i) \geq v(z),
\]
because $f_i$ is lower semi-continuous.
\end{myproof}

\begin{myproof}[Proof of Theorem~\ref{thm:dualitygap}]
Since \eqref{eq:primal} is feasible, $v(0)=p<+\infty$. Let
\[
\xi=\inf\sum_{i=1}^n f_i(x^i),
\]
then by our assumption of $f_i$, $\xi$ is finite. $v(z) \geq \xi$ for all $z \in \RR^m$.  As a consequence, $v(z)$ is bounded below by some affine function, so
\[
-\infty<v^{**}(0) \leq v(0)<+\infty, \quad \epi v^{**}=\cl\conv\epi v.
\]

By Lemma~\ref{lem:perturbfunc}, $v$ is lower semi-continuous. Since $(0,v^{**}(0)) \in \epi v^{**}=\cl\conv\epi v$, for every $\epsilon>0$, there exists $(\kappa,\eta) \in \conv\epi v$ which is sufficiently close to $(0,v^{**}(0))$ such that
\begin{equation}\label{eq:approx}
v(\kappa) \geq v(0)-\epsilon, \quad \eta \leq v^{**}(0)+\epsilon.
\end{equation}
Because $(\kappa,\eta) \in \conv\epi v$, there exists some integer $l$ and $\alpha_j \geq 0$ for $j=1,\dots,l$ such that $\sum_{j=1}^l \alpha_j=1$ and
\[
\kappa=\sum_{j=1}^l\alpha_jz^j, \quad \eta \geq \sum_{j=1}^l\alpha_jv(z^j)
\]
in which $z^j \in \dom v$.

For each $j=1,\dots,l$, find $(\hat x^{1j},\dots,\hat x^{nj})$ attaining the optimal value of the perturbed problem related to $v(z^j)$, i.e.,
\[
v(z^j)=\sum_{i=1}^n f_i(\hat x^{ij}), \quad \sum_{i=1}^n A_i\hat x^{ij} \leq b+z^j,
\]
which means there exists some vector $w^j \in \RR^m_+$ such that
\[
(b+z^j-w^j,v(z^j)) \in \sum_{i=1}^n C_i,
\]
where
\[
C_i=\{(A_ix^i,f_i(x^i)) | f_i(x^i)<+\infty, x^i \in \RR^{n_i}\}.
\]
Taking convex combination of the points above, we have
\[
\left(b+\kappa-\sum_{j=1}^l\alpha_jw^j,\sum_{j=1}^l\alpha_jv(z^j)\right) \in \conv\sum_{i=1}^n C_i.
\]

Now we can apply Corollary~\ref{cor:refinedhalf}, which gives points $(r^i,s_i) \in \conv_{k_i}C_i$ with $1 \leq k_i \leq m+1$ such that
\[
b+\kappa \geq b+\kappa-\sum_{j=1}^l\alpha_jw^j=\sum_{i=1}^n r^i, \quad \eta \geq \sum_{j=1}^l\alpha_jv(z^j) \geq \sum_{i=1}^n s_i
\]
and $\sum_{i=1}^n k_i \leq m+n$. Since $(r^i,s_i) \in \conv_{k_i}C_i$, there exists $\tilde x^{ij} \in \RR^{n_i}$, $\beta_{ij} \geq 0$ for $j=1,\dots,k_i$ such that $f_i(\tilde x^{ij})<+\infty$, $\sum_{j=1}^{k_i} \beta_{ij}=1$ and
\[
r^i=\sum_{j=1}^{k_i} \beta_{ij}A_i\tilde x^{ij}, \quad s_i=\sum_{j=1}^{k_i} \beta_{ij}f_i(\tilde x^{ij}).
\]
Thus,
\begin{equation}\label{eq:feasibilityproof}
\kappa \geq \sum_{i=1}^n\sum_{j=1}^{k_i} \beta_{ij}A_i\tilde x^{ij}-b=\sum_{i=1}^n A_i \sum_{j=1}^{k_i} \beta_{ij}\tilde x^{ij}-b,
\end{equation}
and
\begin{equation}\label{eq:dualitygapproof}
\sum_{i=1}^n \rho^{k_i}_i+\eta \geq \sum_{i=1}^n\left(\rho^{k_i}_i+\sum_{j=1}^{k_i} \beta_{ij}f_i(\tilde x^{ij})\right) \geq \sum_{i=1}^n f_i\left(\sum_{j=1}^{k_i} \beta_{ij}\tilde x^{ij}\right).
\end{equation}

From \eqref{eq:feasibilityproof} we know that
\[
\left(\sum_{j=1}^{k_1} \beta_{1j}\tilde x^{1j},\dots,\sum_{j=1}^{k_n} \beta_{nj}\tilde x^{nj}\right)
\]
is feasible to the perturbed problem related to $v(\kappa)$, so the corresponding objective value
\[
\sum_{i=1}^n f_i\left(\sum_{j=1}^{k_i} \beta_{ij}\tilde x^{ij}\right) \geq v(\kappa).
\]
The above inequality, \eqref{eq:dualitygapproof} and \eqref{eq:approx} imply
\[
v^{**}(0)+\epsilon+\sum_{i=1}^n \rho^{k_i}_i \geq v(0)-\epsilon.
\]

We finish the proof by letting $\epsilon \to 0$ and choose the worst case of $\sum_{i=1}^n \rho^{k_i}_i$ encountered in this process.
\end{myproof}

From a computational viewpoint, since we do not know the $k_i$ that appeared in Theorem~\ref{thm:dualitygap}, in order to find a number for the bound, we have to find the worst case $k_i$ by solving the following optimization problem
\begin{equation}\label{eq:findbound}
\begin{split}
\max \quad & \sum_{i=1}^n \rho^{k_i}_i \\
\st \quad  & 1 \leq k_i \leq m+1, \: k_i \in \ZZ, \quad \forall i=1,\dots,n, \\
           & \sum_{i=1}^n k_i \leq m+n.
\end{split}
\end{equation}

Let $B$ be the optimal value of \eqref{eq:findbound}, then
\[
B \leq \sum_{i=1}^n \rho(f_i).
\]
On the other hand, since for any feasible solution of \eqref{eq:findbound}, the number of $k_i$ with $k_i \geq 2$ is bounded by $m$, so
\[
B=\sum_{i:k_i \geq 2} \rho^{k_i}_i \leq \sum_{i=1}^m \rho(f_i)
\]
if $\rho(f_1) \geq \dots \geq \rho(f_n)$. The above argument shows that the bound $B$ given by the optimization problem \eqref{eq:findbound} is at least as tight as the bound \eqref{eq:udellbound} in \cite{UB2016}.

To illustrate the procedure to calculate the bound $B$, consider the simple case where all the $x^i$ in the primal problem \eqref{eq:primal} are one-dimensional and all the functions $f_i$ equal to the same function $f$. In this case, $\rho^{k_i}_i=\rho(f)$ if $k_i \geq 2$. The optimal value to \eqref{eq:findbound} is attained when the number of $k_i$ that equal to 2 is maximized, so the optimal value is $\min\{m,n\}\rho(f)$, which is the same as the result given by \eqref{eq:udellbound}. The Example~1 used in \cite{UB2016} belongs to this category. It hence explains why the bound \eqref{eq:udellbound} is tight for that example. However, if the dimension of $x^i$ in the primal problem can be arbitrarily large, the bound \eqref{eq:udellbound} can be very loose. As will be shown in Section~\ref{sec:application}, the difference between the bound \eqref{eq:udellbound} and the exact duality gap tends to infinity for a series of problems.

\section{Applications}\label{sec:application}

\subsection{Joint Routing and Congestion Control in Networking}

In this part, we will first apply the previous result to the network utility maximization problem. Consider a network with $N$ users and $L$ links. Let a strictly positive vector $c \in \RR^L$ contain the capacity of each link. Each user $i$ has $K^i$ available paths to send its commodity. We assume that the users are sorted such that $K^1 \geq \dots \geq K^N$. The routing matrix of user $i$, denoted by $R^i$, is a $L \times K^i$ matrix defined by
\[
R^i_{lk}=\begin{dcases*}
1, & if the $k$th path of user $i$ passes through link $l$, \\
0, & otherwise.
\end{dcases*}
\]

Let $x^i \in \RR^{K^i}$ be the vector in which $x^i_k$ is the amount of commodity sent by user $i$ on its $k$th path. Assume that each user $i$ has a utility function $U_i(\cdot)$ depending on the vector $x^i$, then the network utility maximization problem can be written as
\begin{equation}\label{eq:utilitymax}
\begin{split}
\max \quad & \sum_{i=1}^N U_i(x^i) \\
\st \quad  & \sum_{i=1}^N R^ix^i \leq c, \\
           & x^i \geq 0, \quad \forall i=1,\dots,N.
\end{split}
\end{equation}

If all the utility functions $U^i(\cdot)$ are concave, then the above problem \eqref{eq:utilitymax} can be solved by standard convex optimization techniques. Difficulty arises when $U^i(\cdot)$ is not concave. For example, if we restrict each user to choose only one path (single-path routing) and want to maximize the total throughput of the network, then the corresponding utility function is
\[
U_i(x^i)=\max_{s=1,\dots,K^i}x^i_s.
\]
Define
\[
f_i(x^i)=\begin{dcases*}
\min_{s=1,\dots,K^i}(-x^i_s), & if $0 \leq x^i \leq \norm{c}_\infty$, \\
+\infty,                      & otherwise.
\end{dcases*}
\]
Here $\norm{c}_\infty$ is the maximum link capacity in the network. Now the original network utility maximization problem \eqref{eq:utilitymax} is equivalent to the following problem:
\begin{equation}\label{eq:utilitymaxalt}
\begin{split}
\min \quad & \sum_{i=1}^Nf_i(x^i) \\
\st \quad  & \sum_{i=1}^NR^ix^i \leq c.
\end{split}
\end{equation}
The above problem is a particular case of the general optimization problem with separable objectives \eqref{eq:primal} studied in this paper. Using the same technique as shown in Example~\ref{exmp:linear}, we can prove that
\[
\rho^k(f_i) \leq \frac{k-1}{k}\norm{c}_\infty, \quad \rho(f_i)=\frac{K^i-1}{K^i}\norm{c}_\infty.
\]

In the following, suppose each user has a large number of paths to select. More explicitly, $K^i \geq L+1$ is assumed for user $i$. Based on the bound \eqref{eq:udellbound}, the duality gap is bounded by
\[
\sum_{i=1}^{\min\{N,L\}}\frac{K^i-1}{K^i}\norm{c}_\infty,
\]
which is at least
\[
\min\{N,L\}\frac{L}{L+1}\norm{c}_\infty.
\]
In contrast, by Theorem~\ref{thm:dualitygap}, the duality gap is bounded by the optimal value of the following optimization problem:
\[
\begin{split}
\max \quad & \sum_{i=1}^N \frac{k_i-1}{k_i}\norm{c}_\infty \\
\st \quad  & 1 \leq k_i \leq L+1, \: k_i \in \ZZ, \quad \forall i=1,\dots,N, \\
           & \sum_{i=1}^N k_i \leq N+L.
\end{split}
\]
Let $N'$ be the number of users whose $k_i \geq 2$, then $0 \leq N' \leq \min\{N,L\}$. If $N'>0$, using the inequality between arithmetic mean and harmonic mean,
\[
\begin{split}
\sum_{i=1}^N \frac{k_i-1}{k_i}&=\sum_{i:k_i \geq 2} \frac{k_i-1}{k_i}=N'-\sum_{i:k_i \geq 2} \frac{1}{k_i} \\
&\leq N'-\frac{N'^2}{\sum_{i:k_i \geq 2} k_i} \leq N'-\frac{N'^2}{N'+L} \\
&=\frac{L}{1+L/N'} \leq \min\{N,L\}\frac{L}{L+\min\{N,L\}}.
\end{split}
\]
Taking the $N \geq L$ case as an example, by the above inequality, we can bound the duality gap by $L\norm{c}_\infty/2$, essentially half of the bound given by \eqref{eq:udellbound}. The same result was obtained by a specialized technique in \cite{BTT2016}.

Next, we consider another case in which each user has logarithmic utility but still must choose only one path. The utility function of user $i$ can be written as
\[
U_i(x^i)=\log\max_{s=1,\dots,K^i}x^i_s.
\]
Define
\[
g_i(x^i)=\begin{dcases*}
-\log\max_{s=1,\dots,K^i}x^i_s, & if $0 \leq x^i \leq \norm{c}_\infty$, $x^i \neq 0$, \\
+\infty,                       & otherwise.
\end{dcases*}
\]
Then the network utility maximization problem \eqref{eq:utilitymax} is equivalent to the problem obtained by replacing $f_i$ with $g_i$ in \eqref{eq:utilitymaxalt}. Using the result in Example~\ref{exmp:log},
\[
\rho^k(g_i) \leq \log k, \quad \rho(g_i)=\log K^i.
\]

Applying the bound \eqref{eq:udellbound} to this case, we can bound the duality gap by
\begin{equation}\label{eq:oldboundlog}
\sum_{i=1}^{\min\{N,L\}} \log K^i,
\end{equation}
which is at least $\min\{N,L\}\log(L+1)$. On the other hand, by Theorem~\ref{thm:dualitygap}, the duality gap is bounded by the optimal value of the following optimization problem
\begin{equation}\label{eq:loggapmax}
\begin{split}
\max \quad & \sum_{i=1}^N \log{k_i} \\
\st \quad  & 1 \leq k_i \leq L+1, \: k_i \in \ZZ, \quad \forall i=1,\dots,N, \\
           & \sum_{i=1}^N k_i \leq N+L.
\end{split}
\end{equation}
If we still let $N'$ be the number of users whose $k_i \geq 2$, then $0 \leq N' \leq \min\{N,L\}$ and the above bound
\[
\begin{split}
\sum_{i=1}^N \log{k_i}&=\sum_{i:k_i \geq 2} \log{k_i}=\log\prod_{i:k_i \geq 2} k_i \\
                     &\leq \log\left(\frac{\sum_{i:k_i \geq 2} k_i}{N'}\right)^{N'} \leq \log\left(\frac{N'+L}{N'}\right)^{N'} \leq \min\{N,L\}\log\left(1+\frac{L}{\min\{N,L\}}\right),
\end{split}
\]
where in the last step the monotonicity of the function $(1+1/x)^x$ is used. Note that the new bound is qualitatively tighter than the bound \eqref{eq:oldboundlog} provided by \eqref{eq:udellbound}.

\subsection{Dynamic Spectrum Management in Communication}

Consider a communication system consisting of $L$ users sharing a common band. The band is divided equally into $N$ tones. Each user $l$ has a power budget $p_l$ which can be allocated across all the tones. Let $x^i_l$ be the power of user $l$ allocated on tone $i$. Due to the crosstalk interference between users, the total noise for a user on tone $i$ is the sum of a background noise $\sigma_i$ and the power of all other users on the same tone. Therefore, the achievable transmission rate of user $l$ on tone $i$ is given by
\[
u^i_l=\frac{1}{N}\log\left(1+\frac{x^i_l}{\norm{x^i}_1-x^i_l+\sigma_i}\right).
\]

The dynamic spectrum management problem is to maximize the total throughput of all users under the power budget constraints, which can be formulated as the following nonconcave optimization problem:
\begin{equation}\label{eq:dsm}
\begin{split}
\max \quad & \sum_{l=1}^L\sum_{i=1}^Nu^i_l \\
\st \quad  & \sum_{i=1}^N x^i_l \leq p_l, \quad \forall l=1,\dots,L, \\
           & x^i_l \geq 0, \quad \forall i=1,\dots,N,\;\forall l=1,\dots,L.
\end{split}
\end{equation}
For simplicity, we assume that the noises $\sigma_i \leq 1$ and the power budgets $p_l \leq 1$ (if not, then scale all the $\sigma_i$ and $p_l$ simultaneously). The latter requires all the variables $x^i_l \leq 1$. Using the function $h_\sigma$ introduced in Example~\ref{exmp:capacity}, the objective function of \eqref{eq:dsm} can be rewritten as a sum of separable objectives:
\[
\sum_{l=1}^L\sum_{i=1}^Nu^i_l=-\frac{1}{N}\sum_{i=1}^Nh_{\sigma_i}(x^i).
\]

For the purpose of designing dual algorithms, it is of great interest to estimate the duality gap for the problem \eqref{eq:dsm}. In \cite{YL2006}, the authors showed that the duality gap will tend to zero if the number of users $L$ is fixed and the number of tones $N$ goes to infinity. \cite{LZ2009} further determined the convergence rate of the duality gap to be $O(1/\sqrt N)$. Using the bound \eqref{eq:udellbound}, we now demonstrate how to improve the convergence rate estimation to $O(1/N)$, which can be only achieved by the method in \cite{LZ2009} in the special case where all the noises $\sigma_i$ are the same.\footnote{The paper \cite{LZ2009} actually studied the generalization of problem \eqref{eq:dsm} under the existence of path loss coefficient between different users. However, the argument for $O(1/N)$ provided here can also be adapted to the general problem.}

Example~\ref{exmp:capacity} proves that the nonconvexity
\[
\rho^k(h_{\sigma_i}) \leq \log\frac{k}{\sigma_i} \leq \log\frac{k}{\sigma}, \quad \rho(h_{\sigma_i})=\rho^{L+1}(h_{\sigma_i}) \leq \log\frac{L+1}{\sigma},
\]
where $\sigma$ is the minimum among all the noises $\sigma_i$, so \eqref{eq:udellbound} implies that the duality gap is upper bounded by
\begin{equation}\label{eq:oldboundcapacity}
\frac{\min\{N,L\}}{N}\log\frac{L+1}{\sigma},
\end{equation}
which is in the order of $O(1/N)$ if $L$ is fixed and $N$ increases.

In order to further improve the estimation \eqref{eq:oldboundcapacity} for the duality gap, we can resort to Theorem~\ref{thm:dualitygap} and follow the exact same steps for solving \eqref{eq:loggapmax}, which shows that the duality gap is upper bounded by
\[
\frac{\min\{N,L\}}{N}\log\frac{1+L/\min\{N,L\}}{\sigma}.
\]
Like the previous example, our bound is still tighter than the one \eqref{eq:oldboundcapacity} from \eqref{eq:udellbound}.

\section{Generalization with Nonlinear Constraints}\label{sec:nonlinear}

The idea in this paper can also be applied to separable problems with nonlinear constraints such as
\begin{equation}\label{eq:primalnonlinear}
\begin{split}
\min \quad & \sum_{i=1}^n f_i(x^i) \\
\st \quad  & \sum_{i=1}^n g_i(x^i) \leq b.
\end{split}
\end{equation}
Here each $g_i: \RR^{n_i} \to \RR^m$ is a lower semi-continuous function. Note that the problem \eqref{eq:primal} we studied above is a special case of the above optimization problem \eqref{eq:primalnonlinear} if we choose $g_i(x^i)=A_ix^i$. Let $y \in \RR^m_+$ be the dual variables, then the Lagrangian is
\[
L(x,y)=\sum_{i=1}^n (f_i(x^i)+y^Tg_i(x^i))-y^Tb
\]
and the Lagrange dual problem of \eqref{eq:primalnonlinear} is
\[
d=\sup_{y \geq 0}\:\inf_x L(x,y).
\]

If the functions $g_i$ are not convex, the duality gap should not only depend on the nonconvexity of functions $f_i$ but also somehow relate to the functions $g_i$. Like \cite{BS1982}, we want to define the $k$th order nonconvexity of a proper function $f:\RR^n \to \bar\RR$ with respect to another function $g: \RR^n \to \RR^m$, denoted by $\rho^k(f,g)$. To do this, we introduce the following auxiliary function
\[\label{eq:auxfunc}
h^k(x)=\inf_{z \in \RR^n} \left\{f(z) \middle| g(z) \leq \sum_{j=1}^k \beta_jg(x^j), \: \forall \beta_j \geq 0, x^j \in \RR^n \text{ s.t. } \sum_{j=1}^k \beta_j=1 \text{ and } x=\sum_{j=1}^k \beta_jx^j \right\}.
\]
Then $\rho^k(f,g)$ is defined by
\[
\rho^k(f,g)=\sup\left\{h^k\left(\sum_{j=1}^k \alpha_jx^j\right)-\sum_{j=1}^k \alpha_jf(x^j)\right\}.
\]
over all possible convex combinations $\alpha_j \geq 0$, $j=1,\dots,k$, with $\sum_{j=1}^k \alpha_j=1$ of points $x^j$ satisfying $f(x^j)<+\infty$. If the function $g$ is convex, then in the infimum of \eqref{eq:auxfunc} we can choose $z=x$, which gives $h^k(x) \leq f(x)$ and $\rho^k(f,g) \leq \rho^k(f)$. However, the last inequality may not hold if $g$ is not convex.

The proof of Theorem~\ref{thm:dualitygap} can be modified accordingly to the case with nonlinear constraints by replacing $\rho^{k_i}(f_i)$ with $\rho^{k_i}(f_i,g_i)$. In the case when all $g_i$ are convex, $\rho^{k_i}(f_i,g_i) \leq \rho^{k_i}(f_i)$, which implies that the original conclusion in Theorem~\ref{thm:dualitygap} remains true even for convex but nonlinear constraints.

\section{Conclusion}

The improvements obtained in this paper are attributed to two sources. First, instead of using a single number measurement, a series of numbers are introduced to characterize the nonconvexity of a function in a potentially much finer manner. This is based on the concept of $k$th convex hull of a set, which allows us to differentiate different levels of nonconvexity for nonconvex sets. Second, for a separable nonconvex problem, we do not approximate each subproblem individually as people had done before. Instead, by considering all subproblems jointly and noticing that the total deviation of every subproblem to a convex problem is bounded, we reach a much tighter estimation.

\bibliography{SeparableObj}

\begin{thebibliography}{10}
\providecommand{\url}[1]{#1}
\csname url@samestyle\endcsname
\providecommand{\newblock}{\relax}
\providecommand{\bibinfo}[2]{#2}
\providecommand{\BIBentrySTDinterwordspacing}{\spaceskip=0pt\relax}
\providecommand{\BIBentryALTinterwordstretchfactor}{4}
\providecommand{\BIBentryALTinterwordspacing}{\spaceskip=\fontdimen2\font plus
\BIBentryALTinterwordstretchfactor\fontdimen3\font minus
  \fontdimen4\font\relax}
\providecommand{\BIBforeignlanguage}[2]{{%
\expandafter\ifx\csname l@#1\endcsname\relax
\typeout{** WARNING: IEEEtranS.bst: No hyphenation pattern has been}%
\typeout{** loaded for the language `#1'. Using the pattern for}%
\typeout{** the default language instead.}%
\else
\language=\csname l@#1\endcsname
\fi
#2}}
\providecommand{\BIBdecl}{\relax}
\BIBdecl

\bibitem{Asplund1963}
E.~Asplund, ``A {k}-extreme point is the limit of {k}-exposed points,''
  \emph{Israel Journal of Mathematics}, vol.~1, no.~3, pp. 161--162, Sep. 1963.

\bibitem{AE1976}
J.~P. Aubin and I.~Ekeland, ``Estimates of the duality gap in nonconvex
  optimization,'' \emph{Mathematics of Operations Research}, vol.~1, no.~3, pp.
  225--245, Aug. 1976.

\bibitem{BS1982}
D.~P. Bertsekas and N.~R. Sandell, ``Estimates of the duality gap for
  large-scale separable nonconvex optimization problems,'' in \emph{Decision
  and Control, 1982 21st IEEE Conference on}, Dec. 1982, pp. 782--785.

\bibitem{BTT2016}
Y.~Bi, C.~W. Tan, and A.~Tang, ``Network utility maximization with path
  cardinality constraints,'' in \emph{IEEE INFOCOM 2016}, Apr. 2016.

\bibitem{ET1999}
I.~Ekeland and R.~Temam, \emph{Convex Analysis and Variational Problems}, 1999.

\bibitem{HL1993}
J.-B. Hiriart-Urruty and C.~Lemar{\'e}chal, \emph{Convex Analysis and
  Minimization Algorithms II}.\hskip 1em plus 0.5em minus 0.4em\relax Springer
  Berlin Heidelberg, 1993.

\bibitem{LS2013}
J.~Lawrence and V.~Soltan, ``{Carath\'{e}odory}-type results for the sums and
  unions of convex sets,'' \emph{Rocky Mountain Journal of Mathematics},
  vol.~43, no.~5, pp. 1675--1688, Oct. 2013.

\bibitem{LZ2009}
Z.-Q. Luo and S.~Zhang, ``Duality gap estimation and polynomial time
  approximation for optimal spectrum management,'' \emph{IEEE Transactions on
  Signal Processing}, vol.~57, no.~7, pp. 2675--2689, Jul. 2009.

\bibitem{Starr1969}
R.~M. Starr, ``Quasi-equilibria in markets with non-convex preferences,''
  \emph{Econometrica}, vol.~37, no.~1, pp. 25--38, Jan. 1969.

\bibitem{UB2016}
M.~Udell and S.~Boyd, ``Bounding duality gap for separable problems with linear
  constraints,'' \emph{Computational Optimization and Applications}, vol.~64,
  no.~2, pp. 355--378, Jun. 2016.

\bibitem{YL2006}
W.~Yu and R.~Lui, ``Dual methods for nonconvex spectrum optimization of
  multicarrier systems,'' \emph{IEEE Transactions on Communications}, vol.~54,
  no.~7, pp. 1310--1322, Jul. 2006.

\end{thebibliography}
\end{document}